\documentclass[pdflatex,sn-mathphys-ay]{sn-jnl}

\usepackage{graphicx}
\usepackage{multirow}
\usepackage{amsmath,amssymb,amsfonts}
\usepackage{amsthm}
\usepackage{mathrsfs}
\usepackage[title]{appendix}
\usepackage{xcolor}
\usepackage{textcomp}
\usepackage{booktabs}
\usepackage{algorithm}
\usepackage{algorithmicx}
\usepackage{algpseudocode}
\usepackage{subcaption}
\usepackage{listings}
\usepackage{url}
\usepackage{setspace}
\usepackage[acronym]{glossaries}
\usepackage{etoolbox}
\usepackage[authoryear]{natbib}

\makeglossaries

\makeatletter
\patchcmd{\@subsubsec}{\itshape}{}{}{}
\makeatother

\theoremstyle{thmstyleone}

\theoremstyle{thmstyletwo}

\theoremstyle{thmstylethree}

\begin{document}

\newacronym[plural=SCs, firstplural=Supply Chains (SCs)]{sc}{SC}{Supply Chain}

\title[Robust and Adaptive Optimization]{IoT-based Fresh Produce Supply Chain Under Uncertainty: An Adaptive Optimization Framework}

\author*[1]{\fnm{Chirag} \sur{Seth}}\email{cseth@uwaterloo.ca}
\author[1]{\fnm{Mehrdad} \sur{Pirnia}}\email{mpirnia@uwaterloo.ca}
\author[1]{\fnm{James H.} \sur{Bookbinder}}\email{jbookbinder@uwaterloo.ca}

\affil*[1]{\orgdiv{Department of Management Science and Engineering}, \orgname{University of Waterloo}, \orgaddress{\street{200 University Ave W}, \city{Waterloo}, \postcode{N2L 3G1}, \state{Ontario}, \country{Canada}}}

\abstract{
Fruits and vegetables form a vital component of the global economy; however, their distribution poses complex logistical challenges due to high perishability, supply fluctuations, strict quality and safety standards, and environmental sensitivity. In this paper, we propose an adaptive optimization model that accounts for delays, travel time, and associated temperature changes impacting produce shelf life, and compare it against traditional approaches such as Robust Optimization, Distributionally Robust Optimization, and Stochastic Programming. Additionally, we conduct a series of computational experiments using Internet of Things (IoT) sensor data to evaluate the performance of our proposed model. Our study demonstrates that the proposed adaptive model achieves a higher shelf life, extending it by over 18\% compared to traditional optimization models, by dynamically mitigating temperature deviations through a temperature feedback mechanism. The promising results demonstrate the potential of this approach to improve both the freshness and efficiency of logistics systems—an aspect often neglected in previous works.
}


\maketitle

\section{Introduction}

The logistics and management of perishable \glspl{sc} are inherently complex due to factors such as limited shelf-life, demand variability, stringent safety requirements, and sensitivity to environmental conditions. Approximately one-third of global food production, equating to about 1.3 billion tons per year, is lost or wasted along these \glspl{sc} \citep{gustavsson2011globala}. Specifically, losses in fruits and vegetables \glspl{sc} range between 20 to 60 percent globally, influenced by inefficiencies in harvesting, storage, and transportation practices \citep{lemma2014loss}. Developing countries, in particular, face amplified challenges due to inadequate infrastructure, limited cold chain facilities, and lack of modern harvesting techniques, intensifying the complexity of \gls{sc} operations \citep{balaji2016modeling}.  The perishable nature of fresh foods exacerbates their susceptibility to spoilage, demanding precise control of environmental parameters such as temperature and humidity throughout the \gls{sc} \citep{balaji2016modeling}. \\

The Internet of Things (IoT) has emerged as a transformative technology in \gls{sc} management, particularly in managing perishable goods, due to its capability to provide real-time monitoring and decision-making insights. IoT devices such as sensors and actuators, allows continuous tracking of critical parameters like temperature, humidity, and geographic position, crucial for maintaining the quality and safety of perishable products \citep{taj2023iot}. Recent advances in IoT, including sensor-driven shelf-life prediction, real-time logistics re-planning, and automated quality control, have demonstrated significant potential to reduce food spoilage and optimize \gls{sc} efficiency \citep{pang2012value}. Moreover, IoT technologies address common challenges faced by traditional \glspl{sc}, such as information asymmetry and fragmented logistics, thereby enhancing traceability, reliability, and ultimately consumer trust \citep{luo2016intelligent}. Hence, implementation of IoT-based systems in the perishable \gls{sc} not only improves operational efficiency but also offers strategic advantages through increased transparency and responsiveness to dynamic market conditions.
Given the inherent complexities and significant economic impacts associated with perishable goods \glspl{sc}, as well as the promising potential of IoT-based technologies to address these challenges, it is evident that a systematic and integrated approach is crucial. While existing literature emphasizes individual  methods and IoT solutions, there remains a clear gap in the development of comprehensive models that combine  optimization techniques with advanced real-time IoT tracking systems. In the next section, we explore various methods previously published, critically evaluate their contributions, and identify opportunities for further advancements. 

\section{Related Work}
Shelf life estimation is critical in the perishable food \gls{sc} due to its direct impact on food quality, safety, and economic efficiency. Traditional models such as the Arrhenius equation \cite{chen2023} and the  Q$_{10}$ \cite{choi2017persimmons} models are commonly employed to describe the effect of temperature on the degradation rate of fresh produce. The Arrhenius model assumes that reaction rates accelerate exponentially with temperature, this model requires precise estimation of activation energy and reference rates---parameters that are often difficult to determine for a wide variety of fruits and vegetables in heterogeneous \glspl{sc} \cite{sousa2016emerging}. The Q$_{10}$ model, a simplified exponential expression based on temperature differential effects, is thus often used in practice due to its lower data requirements, despite being less physiologically grounded \citep{lamberty2022ambient}. However, both models typically rely on assumptions of uniform ambient conditions and ignore real-time variations that significantly affect spoilage rates. Ambient parameters such as temperature, relative humidity, gas composition (O$_2$ and CO$_2$), and vibration are all known to influence respiration rate, microbial activity, and nutrient degradation in fruits and vegetables \citep{pal2020smart}.\\

Various studies advocate the integration of IoT-enabled sensor systems to track the preceding parameters dynamically across the \gls{sc}. For example, the real-time temperature gradients within the pallets can vary by more than 10\textdegree C, affecting the predicted shelf life by more than 15\% depending on the package position  \citep{pal2020smart}. IoT systems can provide fine-grained ambient profiles, which, when integrated with decay models, enable dynamic shelf-life prediction and more informed logistics decisions. Thus, while classical models laid the foundation, their limitations necessitate hybrid approaches combining empirical modeling with sensor-driven contextual data.\\

Recent empirical research further highlight these modeling challenges. \citep{choi2017persimmons} found that although the Arrhenius model can effectively estimate shelf life under controlled environments, its reliance on accurate activation energy  values and consistent temperature data makes it impractical for variable cold chains \citep{choi2017persimmons}. Their Q$_{10}$ values for semi-dried persimmons varied drastically across temperature bands, reinforcing the need for dynamic adaptation. In addition, those authors emphasized that real-world cold chains suffer from temperature abuse and inconsistencies, which static models cannot account for, thus underlining the significance of integrating IoT for granular monitoring and traceability \citep{aung2023cold}. These findings reinforce that robust shelf-life prediction must account for dynamic, sensor-driven environmental variability across the logistics chain.\\

Optimization methods can play a central role in mitigating spoilage and inefficiencies across perishable \glspl{sc}. Traditional deterministic approaches often fail under real-world variability, prompting researchers to adopt stochastic models that incorporate uncertain demands and travel times \citep{dabbene2008optimisation}. However, stochastic models typically assume known probability distributions, which may not reflect reality, hence the need for robust optimization (RO) and distributionally robust optimization (DRO) approaches. These models protect against worst-case or ambiguous distribution scenarios, this was shown in Govindan et al.\ (2014), who used metaheuristics to solve a multi-objective location-routing problem under sustainability constraints \citep{govindan2014two}. Further, cross-docking-based vehicle scheduling models have been proposed to minimize early/tardy delivery penalties in time-sensitive systems \citep{agustina2014vehicle}, but these often address single-stage or myopic decisions. Our work advances the state of the art by developing an adaptive Optimization framework, enabling dynamic reoptimization as new IoT sensor data becomes available-this bridges a critical gap in temporal decision-making and practical deployability.\\

Other publications further emphasize the growing sophistication of optimization techniques in perishable \glspl{sc}. \citet{krishnan2022robust} proposed a robust multi-objective MILP model for sustainable food networks incorporating perishability, food waste valorization, and supply uncertainty; however, their approach remains limited to static decisions and lacks adaptive responsiveness to real-time data. \citet{qiu2019optimal} developed an exact branch-and-cut algorithm for a perishable production-inventory-routing problem, demonstrating high precision but facing scalability limits and no integration with IoT-based freshness tracking. \citet{hiassat2017genetic} presented a genetic algorithm for the integrated location-inventory-routing problem with perishables, yet did not address re-optimization based on sensor feedback, highlighting a need for dynamic, data-driven models.\citet{musavi2017multi} developed a multi-objective hub location and scheduling model that integrates perishability constraints with sustainability goals by minimizing transportation costs, carbon emissions, and loss of freshness. Their use of a specific metaheuristic provides effective Pareto solutions, though the model assumes fixed fleet capacity and lacks real-time reallocation based on fluctuating conditions.\\

\citet{allaoui2018sustainable} proposed a hybrid approach to sustainable agro-food network design, incorporating economic, environmental, and social criteria. Similarly, \citet{perealopez2003mpc} introduced a rolling-horizon framework for centralized \gls{sc} optimization, which offers a solid foundation for dynamic planning, but does not incorporate IoT feedback loops. \citet{tsang2018intelligent} developed an intelligent IoT-based route planning system that integrates passive cold chain packaging, real-time sensing, and genetic algorithms to minimize food spoilage and optimize routing in a multi-temperature distribution network. While the approach effectively leverages experimental data and Taguchi design, its optimization model operates in a single-stage framework and lacks dynamic feedback adjustment. \citet{peng2019iga} proposed an improved genetic algorithm tailored for cold chain routing under capacity and time window constraints, highlighting significant gains in cost and energy efficiency. However, the model assumes deterministic input and does not incorporate perishability decay or sensor-based variability. Recent advancements further underscore the potential of IoT and advanced optimization in perishable \glspl{sc}. \\

\cite{Shahabi2024} proposed resilience-based models to enhance supply chain flexibility, addressing perishability through sustainable design. \citet{Avishan2025} developed an adaptive optimization approach for production and distribution planning under demand uncertainty, offering dynamic adjustments for perishable food industries. \citet{Shadkam2024} reviewed simulation optimization techniques, advocating for agent-based intelligent digital frameworks to handle supply chain complexity. \citet{Javadi2024} demonstrated IoT applications in dairy supply chains, emphasizing uncertainty handling through real-time monitoring. \citet{Daneshvar2023} introduced a robust possibilistic model for agricultural product distribution under uncertainty, enhancing decision-making flexibility. \citet{Abbasian2022} proposed a hybrid optimization method for sustainable and resilient supply chains, integrating perishability and environmental goals. \citet{Wu2023} explored digital twins for real-time optimization, leveraging IoT data for dynamic perishable supply chain management. \citet{Mirabelli2022} reviewed optimization strategies, highlighting integrated management approaches for perishables. \citet{Ietto2024} introduced an adaptive model predictive control policy for resilient supply chains, addressing uncertain forecasts in perishable goods. Finally, \citet{Pal2020} emphasized smart sensing and communication, enabling precise monitoring and control in perishable food logistics. These studies collectively highlight the need for integrated, adaptive, and IoT-driven optimization frameworks, which our work addresses through dynamic reoptimization and sensor-based decay modeling.Therefore, the main contributions of this paper are as follows:
\begin{itemize}
    \item Proposing an adaptive optimization model that dynamically adjusts routing and temperature control decisions based on temperature changes. \\
    
    \item Employing shelf-life models (Arrhenius and Q$_{10}$) that enable realistic and flexible decay estimation, based on inputs from IoT sensors.\\
    
    \item Conducting an extensive computational comparisons using synthetic IoT-driven datasets to benchmark our proposed approach against  optimization approaches: Deterministic, Robust Optimization (RO), Stochastic Programming (SP) and Distributionally Robust Optimization (DRO).\\
\end{itemize}

\section{Methodology}

This section presents the modeling and optimization framework used to address the transportation of perishable food products under uncertainty. We begin by introducing the temperature-sensitive shelf-life estimation models, namely the Arrhenius and Q$_{10}$ models, which allow for the quantification of freshness degradation during transit. Based on these estimates, we formulate and analyze five optimization models: (i)  deterministic model assuming perfect information, (ii) robust optimization model designed to handle worst-case parameter deviations, (iii)  stochastic programming approach leveraging scenario-based uncertainties, (iv) adistributionally robust optimization (DRO) model for ambiguity-aware decision-making, and (v) the proposed adaptive optimization model that dynamically updates routing and temperature decisions based on IoT sensor feedback. \\

As evident by the previous section, the fresh food \gls{sc} requires careful monitoring and adaptive techniques to maintain quality and minimize waste. The perishability of these goods introduces complexities such as temperature fluctuations, shelf-life degradation, and unpredictable delays. Typically, distributors pre-plan delivery routes based on the expected requirements for the following day. While this approach simplifies logistics, it's rigid route may lead to inefficiencies in terms of both food wastage and total travel time, due to unexpected events.In this paper, perishable products that need to be delivered to various stops while ensuring minimal spoilage and efficient transportation are considered. Table~\ref{tab:variables} includes variables related to temperature dynamics, product perishability, routing decisions, and logistical parameters. For clarity and ease of understanding, each of these notations will be introduced and explained in detail as they appear in the respective sections of the methodology.
\begin{table}[htbp]
\centering
\caption{Key sets, decision variables, and parameters used in shelf-life and transportation modeling.}
\begin{tabular}{lll}
\toprule
\textbf{Symbol} & \textbf{Description} & \textbf{Units} \\
\midrule
\multicolumn{3}{l}{\textit{Sets}} \\
$\mathcal{N}$     & Set of all nodes (warehouse + delivery stops)             & -- \\
$\mathcal{K}$     & Set of perishable product types                           & -- \\
$S$               & Set of scenarios (stochastic model)                       & -- \\
$\Xi$             & Uncertainty set (robust model)                            & -- \\
$\mathcal{P}$     & Ambiguity set of probability distributions (DRO)          & -- \\
\midrule
\multicolumn{3}{l}{\textit{Decision Variables}} \\
$x_{ij}$          & 1 if vehicle travels from node $i$ to $j$, 0 otherwise    & binary \\
$u_i$             & Sequencing variable for node $i$ (MTZ formulation)        & integer \\
$t_{ik}$          & Temperature of product $k$ at node $i$                    & °C \\
$D_{ik}$          & Absolute deviation from ideal temperature $\theta_k$      & °C \\
$S_{ik}$          & Slack variable for temperature bounds                     & °C \\
$I_{tk}$          & Inventory of product $k$ at time $t$                      & units (e.g., kg) \\
\midrule
\multicolumn{3}{l}{\textit{Parameters}} \\
$t$               & Actual temperature                                        & K \\
$t_0$             & Reference temperature                                     & K \\
$L_T$             & Shelf life at temperature $T$                             & hours \\
$L_{T_0}$         & Shelf life at reference temperature $T_0$                 & hours \\
$E_a$             & Activation energy                                         & J/mol \\
$R$               & Universal gas constant ($\approx 8.314$)                  & J/mol$\cdot$K \\
$k$, $k_0$        & Reaction rate constant and pre-exponential factor         & 1/hour \\
$Q_{10}$          &\( Q_{10} \) temperature coefficient                               & dimensionless \\
$T_{ij}$          & Travel time between node $i$ and $j$                      & hours \\
$d_{ij}$          & Distance between node $i$ and $j$                         & km \\
$\delta_i$        & Delay at node $i$                                         & hours \\
$\tau_{ik}$       & Ambient temperature shift at node $i$ for product $k$     & °C \\
$\theta_k$        & Ideal storage temperature for product $k$                 & °C \\
$\theta_k^{\min}, \theta_k^{\max}$ & Min/max allowable temperatures for product $k$   & °C \\
$L_k$             & Initial shelf life of product $k$                         & hours \\
$q_k$             & Demand for product $k$                                    & units (e.g., kg) \\
$\mu_{T_{ij}}, \sigma_{T_{ij}}^2$ & Mean and variance of travel time from $i$ to $j$     & hours, hours$^2$ \\
$\mu_{\delta_i}, \sigma_{\delta_i}^2$ & Mean and variance of delay at node $i$         & hours, hours$^2$ \\
$\beta$           & Temperature correction factor                             & [0,1] \\
$\lambda_1, \lambda_2$ & Penalty weights for temperature deviation and slack     & cost units \\
\bottomrule
\end{tabular}
\label{tab:variables}
\end{table}

\subsection{Shelf Life Calculation}

Since degradation rates depend significantly on temperature and environmental factors, a dynamic shelf-life estimation model would allow real-time decision-making at key points such as delivery stops or after door openings. By incorporating a predictive model, we can ensure that deliveries are prioritized based on the freshness of each item, leading to reduced waste and improved efficiency. To quantify the impact of temperature variations on shelf life, we will use two models: the Arrhenius model and the Q$_{10}$ model. The Arrhenius model \citep{arrhenius1889} is a mechanistic approach that explains how reaction rates, including food spoilage processes, are influenced by temperature. It is based on the Arrhenius equation:

\begin{equation}
k = k_0 \exp \left( -\frac{E_a}{R T} \right),
\label{eq:arrhenius}
\end{equation}

where \( k \) is the temperature-dependent degradation rate constant of the food quality index, \( k_0 \) is the pre-exponential (frequency) factor, \( E_a \) is the activation energy (in J/mol) required to initiate the degradation process, \( R \) is the universal gas constant (\( 8.314 \) J/mol·K), and \( T \) is the absolute temperature (in Kelvin) at which spoilage occurs. \\

Using this relationship, we can estimate the remaining shelf life at a given temperature by applying:

\begin{equation}
L_T = L_{T_0} \exp \left[\frac{E_a}{R} \left( \frac{1}{T} - \frac{1}{T_0} \right) \right],
\label{eq:arrhenius_shelf}
\end{equation}

where \( L_T \) denotes the estimated shelf life of the product at temperature \( T \), and \( L_{T_0} \) is the known shelf life under a reference storage temperature \( T_0 \). This exponential relationship models the accelerated degradation of perishable products under elevated temperatures.\\

The Arrhenius model is highly accurate for predicting spoilage behavior when the chemical kinetics are well characterized for that product. However, its practical application can be challenging as it requires product-specific activation energy values, which may not always be available for a wide range of fresh produce.\\

As an alternative, the Q$_{10}$ model provides a simpler empirical approach to estimating shelf life. It assumes that for every \( 10^\circ\text{C} \) rise in temperature, the rate of deterioration increases by a factor of \( Q_{10} \), typically ranging between 1.5 and 3 for food products. The shelf-life relationship is given by:

\begin{equation}
L_T = L_{T_0} Q_{10}^{\frac{T_0 - T}{10}},
\label{eq:Q$_{10}$_shelf}
\end{equation}

where \( Q_{10} \) is the temperature coefficient for spoilage acceleration, and \( T \), \( T_0 \), \( L_T \), and \( L_{T_0} \) are defined as above. Unlike the Arrhenius model, the Q$_{10}$ model does not require detailed knowledge of activation energies, making it more practical for applications where specific chemical properties of the produce are unknown. One can also modify this to \( Q_5 \) or \( Q_3 \) based on the application domain, especially when the temperature variation is relatively small or constrained. However, the Q$_{10}$ model assumes a fixed exponential relationship between temperature and spoilage, which may not always be accurate for certain produce that exhibit non-linear degradation behavior under fluctuating environmental conditions. In this work, the Arrhenius model is used when sufficient product-specific data is available to provide precise degradation estimates, while the Q$_{10}$ model is employed as a more adaptable alternative when dealing with a variety of produce types where exact kinetic parameters are unknown.\\

\subsection{Optimization Models for Fresh Produce Supply Chains}
 In the following sections, we explore a sequence of mathematical optimization models that are each suited to different operational contexts and data availability scenarios. These models incorporate perishability constraints either explicitly through shelf-life degradation estimates or implicitly via temperature-based thresholds to assess the quality and viability of delivery schedules. Shelf-life calculations derived from the Arrhenius or Q$_{10}$ models serve as critical inputs for quantifying freshness and prioritizing deliveries.

\subsubsection{Deterministic Model}

Assuming perfect knowledge of all parameters, a deterministic model provides a baseline formulation for the transportation planning problem. In this case, the objective is to minimize the total travel cost while satisfying operational and perishability constraints. The generic structure of a deterministic optimization problem is given by:\\

\begin{subequations}
\begin{align}
[\mathcal{D}]: \quad & \min_{x \in X} f(x, \bar{\xi}) \label{eq:det_obj} \\
& \text{s.t.} \quad f_i(x, \bar{\xi}) \leq 0 \quad \forall i \in I, \label{eq:det_const}
\end{align}
\end{subequations}

where \( x \) represents the decision variables (such as routing choices), \( \bar{\xi} \) is the known vector of input parameters (including travel times, delays, and shelf-life values), and \( I \) is the index set of constraints.  Based on this structure, the specific formulation of the transportation problem is as follows:\\

\begin{equation}
\min \sum_{i \in N} \sum_{j \in N} (t_{ij} + \delta_i) x_{ij}
\label{eq:det_min}
\end{equation}

subject to

\begin{align}
\sum_{j \in N \setminus \{0\}} x_{0j} &= 1 \label{eq:det_start} \\
\sum_{i \in N \setminus \{0\}} x_{i0} &= 1 \label{eq:det_end} \\
\sum_{j \in N \setminus \{i\}} x_{ij} &= 1 \quad \forall i \in N \setminus \{0\} \label{eq:det_visit_out} \\
\sum_{j \in N \setminus \{i\}} x_{ji} &= 1 \quad \forall i \in N \setminus \{0\} \label{eq:det_visit_in} \\
x_{ii} &= 0 \quad \forall i \in N \label{eq:det_noloop} \\
u_i - u_j + n x_{ij} &\leq n - 1 \quad \forall i \neq j, \; i, j \in N \setminus \{0\} \label{eq:det_subtour} \\
\sum_{i \in N} \sum_{j \in N} ( \delta_i + t_{ij} ) x_{ij} &\leq L_k - R_k \quad \forall k \in K \label{eq:det_shelf}
\end{align}

\vspace{1em}

The objective function in Equation~\eqref{eq:det_min} minimizes the total transportation cost, which includes the travel time \( t_{ij} \) between nodes \( i \) and \( j \), and the service delay \( \delta_i \) incurred at node \( i \). The binary decision variable \( x_{ij} \in \{0,1\} \) indicates whether the vehicle moves directly from node \( i \) to node \( j \).\\

Constraints~\eqref{eq:det_start} and~\eqref{eq:det_end} ensure that the delivery route originates and terminates at the depot, denoted by node 0. Constraints~\eqref{eq:det_visit_out} and~\eqref{eq:det_visit_in} guarantee that each customer node \( i \in N \setminus \{0\} \) is visited exactly once, with a unique incoming and outgoing arc.\\

Constraint~\eqref{eq:det_noloop} eliminates self-loops by enforcing \( x_{ii} = 0 \) for all \( i \in N \), ensuring that the vehicle does not return to the same node without progressing. Constraint~\eqref{eq:det_subtour} implements the subtour elimination rule using the Miller–Tucker–Zemlin (MTZ) formulation. Here, \( u_i \in \mathbb{R}^{+} \) is a sequencing variable that defines the order of visiting node \( i \), and \( n \) is the total number of customer nodes. Constraint~\eqref{eq:det_shelf} imposes a perishability limit. For each product type \( k \in \mathcal{K} \), the total time spent in transit, computed as the sum of all \( t_{ij} \) and \( \delta_i \) values along the selected route, must not exceed the available shelf life window \( L_k - R_k \), where \( L_k \) is the initial shelf life at dispatch and \( R_k \) is the minimum required remaining shelf life upon delivery, assuming \( L_k - R_k  >0\).\\

This deterministic model serves as a foundational formulation for routing perishable goods under idealized conditions, where all parameters—travel times \( t_{ij} \), delays \( \delta_i \), and shelf-life values \( L_k, R_k \)—are known in advance and remain fixed throughout the planning horizon.

\subsubsection{Robust Optimization Model}

This section introduces the concept of robust optimization, detailing its principles and methodologies \citep{sun2017robust}. Initially proposed by Soyster \citep{soyster1973convex} and further advanced by Ben-Tal and Nemirovski \citep{ben1998robust,ben1999robust,ben2000robust}, robust optimization has become a vital approach for tackling problems under uncertainty. By accounting for uncertainties and variations in input parameters, robust optimization offers a framework for optimization under uncertainty, ensuring system performance and reliability despite unpredictable factors. Consider the general form of a deterministic optimization problem as defined in~\eqref{eq:det_obj}--\eqref{eq:det_const}, where \( x \) represents the decision variables, \( X \subset \mathbb{R}^n \) is the feasible set, \( \bar{\xi} \) is the vector of nominal values of the uncertain parameters \( \xi \), and \( I = \{1, \dots, m\} \) is the index set of constraints. This problem seeks the optimal \( x \) to minimize the objective function \( f(x, \bar{\xi}) \) under constraints \( f_i(x, \bar{\xi}) \leq 0, \; \forall i \in I \).\\

In an optimization problem, some parameters such as travel and handling times are inherently non-deterministic. When \( \xi \neq \bar{\xi} \), the optimal solution obtained may not be feasible for some constraints \( f_i(x, \xi) \leq 0 \), or even if all constraints are satisfied, the objective function value may vary for different \( \xi \). This limitation leads to the adoption of robust optimization, which aims to minimize the worst-case objective function and ensure solution feasibility across all possible realizations of uncertain parameters within a predefined uncertainty set. \\The robust counterpart of the preceding problem is formulated as follows:\\

\begin{subequations}
\begin{align}
[\mathcal{RO}]: \quad & \min_{x \in X} \sup_{\xi \in \Xi} f(x, \xi) \label{eq:ro_obj} \\
& \text{s.t.} \quad f_i(x, \xi) \leq 0 \quad \forall \xi \in \Xi, \; \forall i \in I, \label{eq:ro_const}
\end{align}
\end{subequations}

where \( \Xi \) is the uncertainty set of the parameter \( \xi \). That is, \( \Xi \) represents all possible values that may be assumed by the random parameter \( \xi \). This framework minimizes the worst-case objective in~\eqref{eq:ro_obj} under the constraints in~\eqref{eq:ro_const}, ensuring that a solution \( x \) is feasible for all realizations of the uncertain parameter \( \xi \) within the uncertainty set \( \Xi \). The tractable reformulation of the problem depends on both the functions \( f \) and the structure of the uncertainty set \( \Xi \). We define the robust perishable delivery optimization problem considering the worst-case uncertainties in travel time, delays, and shelf life variations. The goal is to minimize the transportation cost while ensuring that the freshness of the products is preserved during transportation, even under uncertain conditions. The uncertainty set used in this study is a box uncertainty set, where the travel times and delays are bounded within given minimum and maximum values.

The mathematical model is given as:

\begin{equation}
\min \sum_{i \in N} \sum_{j \in N} (T_{\text{max}, ij} + \delta_{\text{max}, i}) x_{ij}
\label{eq:ro_min}
\end{equation}

In Equation~\eqref{eq:ro_min}, \( T_{\text{max},ij} \) denotes the maximum possible travel time between nodes \( i \) and \( j \), and \( \delta_{\text{max},i} \) denotes the maximum possible delay at location \( i \). The variable \( x_{ij} \in \{0,1\} \) is the binary routing decision variable. Subject to constraints~\eqref{eq:det_start}--\eqref{eq:det_subtour}, which enforce routing feasibility, and:

\begin{equation}
\sum_{i \in N} \sum_{j \in N} ( \delta_{\text{max}, i} + T_{\text{max}, ij} ) x_{ij} \leq L_k - R_k \quad \forall k \in K
\label{eq:ro_shelf}
\end{equation}

Constraint~\eqref{eq:ro_shelf} ensures that, for each product \( k \in \mathcal{K} \), the total transit time in the worst-case scenario (including worst-case delays \( \delta_{\text{max},i} \) and travel times \( T_{\text{max},ij} \)) does not exceed the freshness window \( L_k - R_k \) The uncertainty set for travel times and delays is defined as:

\begin{subequations}
\begin{align}
T_{ij} \in [ T_{ij}^{\min}, T_{ij}^{\max} ], \quad \forall i \in N, \forall j \in N \label{eq:travel_uncertainty} \\
\delta_i \in [ \delta_i^{\min}, \delta_i^{\max} ], \quad \forall i \in N \label{eq:delay_uncertainty}
\end{align}
\end{subequations}

Here, \( T_{ij}^{\min} \) and \( T_{ij}^{\max} \) represent the minimum and maximum possible travel times between nodes \( i \) and \( j \), respectively. Similarly, \( \delta_i^{\min} \) and \( \delta_i^{\max} \) define the minimum and maximum possible delays at location \( i \).

\subsubsection{Stochastic Programming Model}

Real-world transportation involves uncertainty in travel times, demand, and handling delays. An alternative method to handle randomness is Stochastic Programming (SP), introduced by \citet{dantzig1955linear}. SP considers discrete scenarios for each parameter and optimizes the expected value based on their probability distributions, derived from historical data, leading to the following formulation:

\begin{equation}
\min_{x} \mathbb{E}_{\xi} \left[ f(x, \xi) \right]
\label{eq:sp_obj}
\end{equation}

In this formulation, \( x \) represents the vector of decision variables, \( \xi \) is the random vector representing uncertainty in system parameters, and \( f(x, \xi) \) is the cost function that depends on both decision variables and realizations of uncertainty. The operator \( \mathbb{E}_{\xi} \left[ \cdot \right] \) denotes the expectation over all realizations of \( \xi \), thereby accounting for the expected transportation cost across possible scenarios.\\

In Stochastic Programming, uncertainty in travel times \( t_{ij} \) and delays \( \delta_i \) is modeled by a finite set of discrete scenarios. Each scenario \( s \in S \) represents a particular realization of these uncertain parameters, with associated probability \( p_s \in [0,1] \), where \( \sum_{s \in S} p_s = 1 \). The travel time and delay values for scenario \( s \) are denoted by \( t_{ij}^s \) and \( \delta_i^s \), respectively. The model optimizes the expected value of the total transportation cost while respecting perishability constraints. The base model uses the objective function given by~\eqref{eq:det_min} and routing constraints~\eqref{eq:det_start}--\eqref{eq:det_subtour}, but replaces the in-transit shelf-life constraint with the following expected constraint:

\begin{equation}
\sum_{s \in S} p_s \sum_{i \in N} \sum_{j \in N} \left( \delta_i^s + t_{ij}^s \right) x_{ij} \leq L_k - R_k \quad \forall k \in K
\label{eq:sp_shelf}
\end{equation}

In this constraint, \( L_k \geq 0 \) is the initial shelf life available for product \( k \), and \( R_k \geq 0 \) is the minimum required shelf life upon delivery. The binary variable \( x_{ij} \in \{0,1\} \) indicates whether the route from node \( i \) to node \( j \) is selected. The constraint in Equation~\eqref{eq:sp_shelf} ensures that, on average across all defined scenarios \( s \in S \), the total time spent in transit—including scenario-specific delays \( \delta_i^s \) and travel times \( t_{ij}^s \)—does not reduce the usable shelf life of any product \( k \) below the required level. This formulation balances expected performance with perishability constraints by incorporating probabilistic modeling of operational uncertainty. The risk aversion of the solution can be tuned by introducing a robustness parameter \( z \) see~\eqref{eq:dro_shelf}, which can modify the constraint threshold or enter the objective in extended formulations. This approach optimizes the transportation process while accounting for uncertainty in travel times and delays, ensuring product freshness is preserved.\\

\subsubsection{Distributionally Robust Model}

Distributionally Robust Optimization (DRO) extends stochastic programming by considering a set of probability distributions rather than a single fixed distribution \citep{bertsimas2019adaptive}. The goal is to minimize the worst-case expected cost:

\begin{equation}
\min_{x} \sup_{\mathbb{P} \in \mathcal{P}} \mathbb{E}_{\mathbb{P}} \left[ f(x, \xi) \right],
\label{eq:dro_obj}
\end{equation}

where \( \mathcal{P} \) is an ambiguity set of probability distributions. DRO is useful in \gls{sc} logistics where disruptions such as weather, strikes, or demand surges lead to distributional uncertainty. The tractable reformulation of the problem depends on the structure of the ambiguity set. In the present research, a moment-based ambiguity set is used, which considers the mean and variance of the uncertain parameters, such as travel times \( T_{ij} \) and delays \( \delta_i \). The ambiguity set is defined as follows:

\begin{equation}
\mathcal{P} = \left\{ P \mid \mathbb{E}_P[T_{ij}] = \mu_{T_{ij}}, \mathbb{E}_P[\delta_i] = \mu_{\delta_i}, \text{Var}_P[T_{ij}] = \sigma_{T_{ij}}^2, \text{Var}_P[\delta_i] = \sigma_{\delta_i}^2 \right\},
\label{eq:dro_ambiguity}
\end{equation}

In this ambiguity set, \( \mu_{T_{ij}} \) and \( \sigma_{T_{ij}}^2 \) represent the mean and variance of the travel time between locations \( i \) and \( j \), while \( \mu_{\delta_i} \) and \( \sigma_{\delta_i}^2 \) represent the mean and variance of the delay time at location \( i \). The mathematical model is given by~\eqref{eq:det_min}, subject to constraints~\eqref{eq:det_start}--\eqref{eq:det_subtour} and:

\begin{equation}
\sum_{i \in N} \sum_{j \in N} \left( \mu_{\delta_i} + \mu_{T_{ij}} \right) x_{ij} + z \cdot \sum_{i \in N} \sum_{j \in N} \left( \sigma_{\delta_i}^2 + \sigma_{T_{ij}}^2 \right) x_{ij}^2 \leq Q_k \cdot L_k - R_k \quad \forall k \in K
\label{eq:dro_shelf}
\end{equation}

  In this model, \( x_{ij} \in \{0,1\} \) is the decision variable representing the link selection between locations \( i \) and \( j \), \( \mu_{\delta_i} \) and \( \mu_{T_{ij}} \) are the mean delay times and travel times, \( \sigma_{\delta_i}^2 \) and \( \sigma_{T_{ij}}^2 \) are the variances of the delay times and travel times, \( Q_k \geq 0 \) represents the demand for product \( k \), \( L_k \geq 0 \) represents the available shelf life for product \( k \), \( R_k \geq 0 \) represents the required shelf life for product \( k \), and \( z \geq 0 \) is a risk-aversion parameter or safety factor, which controls the level of robustness against uncertainty in the random parameters. The moment-based DRO in-transit time constraint~\eqref{eq:dro_shelf} ensures that the total time spent in transit does not exceed the available shelf life \( L_k \) minus the required shelf life \( R_k \) for each product, considering both the mean and variance of the travel times and delays. As always, constraints~\eqref{eq:det_start} and~\eqref{eq:det_end} enforce that the route starts and ends at the warehouse.

\subsection{The Proposed Adaptive Optimization Model}

Our adaptive perishable goods transportation model optimizes the delivery of perishable produce across multiple locations while ensuring that temperature and shelf life constraints are met at the time of delivery. The goal is to minimize the total transportation cost while maintaining the freshness of products and adhering to temperature and time constraints.\\

In real-world fresh-food logistics, environmental conditions, such as ambient temperature and handling delays, introduce dynamic and location-specific effects on product quality. Unlike static formulations in Section 3.2, this model incorporates temperature feedback and adjustment mechanisms across sequential stops to minimize spoilage and energy costs while preserving product freshness. The model minimizes total cost considering travel time, temperature deviations, and operational flexibility. In practice, this optimization is carried out iteratively, where at each visited node \(i\), a restricted MILP is solved to determine the next stop \(j\), and the temperature state is updated using feedback from current conditions. At each iteration, we define a dynamic subset \(\textit{Remaining} \subset \mathcal{N} \setminus \textit{Visited}\), which contains the nodes (delivery stops) that are yet to be visited. This set is updated after each optimization step by removing the node selected in the previous hop. The decision model then selects the best next node \(j \in \textit{Remaining}\) to visit from the current node \(i\), minimizing the travel and freshness penalty locally at each stage. The objective and constraints in the following apply to each such sub-problem in the adaptive sequence. The mathematical formulation for each decision step is given as follows:

\begin{equation}
\min_{j \in \textit{Remaining}} \left( T_{ij} + \delta_i \right) x_{ij} + \lambda_1 \sum_{k \in \mathcal{K}} D_{jk} + \lambda_2 \sum_{k \in \mathcal{K}} S_{jk}
\label{eq:ms_min}
\end{equation}

subject to constraints~\eqref{eq:det_start}--\eqref{eq:det_subtour} and:

\begin{align}
t_{0k} &= \theta_k \quad \forall k \in \mathcal{K} \label{eq:ms_temp_init} \\
t_{jk} &\geq t_{ik} + \tau_{jk} + \beta (\theta_k - t_{ik}) - M (1 - x_{ij}) - S_{jk} \quad \forall i, j \in \mathcal{N}, i \neq j, k \in \mathcal{K} \label{eq:ms_temp_lower} \\
t_{jk} &\leq t_{ik} + \tau_{jk} + \beta (\theta_k - t_{ik}) + M (1 - x_{ij}) + S_{jk} \quad \forall i, j \in \mathcal{N}, i \neq j, k \in \mathcal{K} \label{eq:ms_temp_upper} \\
\theta_k^{\min} &\leq t_{ik} \leq \theta_k^{\max} \quad \forall i \in \mathcal{N}, k \in \mathcal{K} \label{eq:ms_temp_bounds} \\
D_{ik} &\geq t_{ik} - \theta_k \quad \forall i \in \mathcal{N}, k \in \mathcal{K} \label{eq:ms_dev_pos} \\
D_{ik} &\geq \theta_k - t_{ik} \quad \forall i \in \mathcal{N}, k \in \mathcal{K} \label{eq:ms_dev_neg}
\end{align}

This model includes a dynamic temperature update mechanism that accounts for ambient shifts \( \tau_{jk} \), moderated by the temperature correction factor \( \beta \), along each traveled arc \( (i, j) \). In this formulation:
\( t_{0k} \) is the initial temperature of product \( k \) at the warehouse (node 0),
\( t_{jk} \) is the temperature of product \( k \) at location \( j \),
\( \tau_{jk} \) denotes the ambient temperature shift along the arc from node \( i \) to \( j \),
\( \beta \in [0,1] \) is the temperature correction factor moderating adjustment toward the required temperature \( \theta_k \),
\( M \) is a large positive constant used for constraint relaxation when the arc \( (i,j) \) is not selected (i.e., when \( x_{ij} = 0 \)),
\( S_{jk} \geq 0 \) is a slack variable allowing temperature deviations to maintain feasibility,
\( \theta_k \) is the required storage temperature for product \( k \),
\( \theta_k^{\min} \) and \( \theta_k^{\max} \) are the minimum and maximum allowable temperatures for product \( k \),
\( D_{ik} \) represents the absolute temperature deviation of product \( k \) at node \( i \),
\( \lambda_1, \lambda_2 \geq 0 \) are penalty coefficients on temperature deviations and slack variables, respectively,
\( x_{ij} \in \{0,1\} \) is the route decision variable indicating if arc \( (i,j) \) is traversed,
\( T_{ij} \) is the travel time between nodes \( i \) and \( j \),
\( \delta_i \) is the delay at node \( i \),
\( \mathcal{N} \) is the set of all nodes (including warehouse and delivery stops),
\( \mathcal{K} \) is the set of all product types.\\

The objective in Equation~\eqref{eq:ms_min} minimizes local travel time and delay costs at each hop, penalizes temperature deviations weighted by \( \lambda_1 \), and penalizes slack usage weighted by \( \lambda_2 \), over all candidate destinations \( j \in \textit{Remaining} \). The total route cost is accumulated incrementally across all stages of the rolling-horizon execution. Routing feasibility is enforced by constraints~\eqref{eq:det_start} to~\eqref{eq:det_subtour}, which respectively ensure route start/end at the warehouse, single visits to each delivery stop, elimination of self-loops, and prevention of subtours via the MTZ formulation.\\

Constraint~\eqref{eq:ms_temp_init} sets the initial temperature of each product at the warehouse. In each subsequent MILP re-optimization, the current temperature \( t_{ik} \) serves as the initial state for that subproblem. Constraints~\eqref{eq:ms_temp_lower} and~\eqref{eq:ms_temp_upper} impose lower and upper bounds on temperature evolution between successive nodes, with the big-\( M \) term deactivating constraints on arcs not selected in the route. Constraint~\eqref{eq:ms_temp_bounds} restricts product temperatures within allowable limits. Constraints~\eqref{eq:ms_dev_pos} and~\eqref{eq:ms_dev_neg} compute the absolute deviation \( D_{ik} \) from the target temperature \( \theta_k \).\\

\subsubsection{Adaptive Routing Procedure}
To better illustrate the workings of the proposed model, we present a stepwise algorithm \autoref{alg:adaptive_routing} that reflects the actual rolling-horizon implementation. At each stage of the delivery process, the model evaluates only the remaining candidate nodes that are yet to be visited. This subset is referred to as \textit{Remaining}, and is updated after each decision.\\

The optimization at each step solves a restricted MILP based on the current node \( i \), and selects the best next node \( j \in \textit{Remaining} \) based on travel time, delay, temperature drift, and expected freshness penalties. The temperature state is propagated forward using the dynamic correction mechanism described in Equations~\eqref{eq:ms_temp_lower}--\eqref{eq:ms_temp_upper}. The process continues until all stops have been visited, and the route concludes with a return to the warehouse.\\

\begingroup
\setstretch{1}  
\begin{algorithm}[H]
\caption{Rolling–Horizon Adaptive Optimization with Temperature Feedback}
\label{alg:adaptive_routing}
\begin{algorithmic}[1]

\Require 
Node set $\mathcal{N} = \{0,1,\dots,n-1\}$ (warehouse $0$ and customers),\\
Product set $\mathcal{K} = \{1,\dots,K\}$, travel times $T_{ij}$, delays $\delta_i$,\\
Ambient temperature shifts $\tau_{jk}$, temperature bounds $\left[\theta_k^{\min}, \theta_k^{\max}\right]$,\\
Ideal temperatures $\theta_k$, correction factor $\beta$, penalties $\lambda_1, \lambda_2$, big-$M$

\Ensure 
Route $\{x_{ij}=1\}$, temperature profile $\{t_{ik}\}$, total cost

\Statex \textbf{Phase 0: Initialization}
\State Set $t_{0k} \gets \theta_k$ \quad for all $k \in \mathcal{K}$ \hfill \emph{(initial temperature)}
\State $\textit{Visited} \gets \{0\}$, \quad $\textit{Remaining} \gets \mathcal{N} \setminus \{0\}$
\State Initialize cost: $\textit{TotalCost} \gets 0$, route: $\textit{Route} \gets \emptyset$

\Statex \textbf{Phase 1–$n$: Iterative Re-Optimization and Feedback}
\While{$\textit{Remaining} \neq \emptyset$}
    \State Let $i \gets$ most recently visited node
    \State Solve MILP defined by constraints~\eqref{eq:det_start}--\eqref{eq:det_subtour}, \eqref{eq:ms_temp_init}--\eqref{eq:ms_dev_neg}, restricted to nodes $\textit{Visited} \cup \textit{Remaining}$

    \State Select next node $j \in \textit{Remaining}$ such that $x_{ij} = 1$

    \ForAll{$k \in \mathcal{K}$}
        \State Compute drift $\tau_{jk}$ and correction $\beta(\theta_k - t_{ik})$
        \State Predict temperature at $j$:
        \[
        t_{jk} \gets t_{ik} + \tau_{jk} + \beta(\theta_k - t_{ik})
        \]
        \State Impose constraints:
        \begin{align*}
        t_{jk} &\ge t_{ik} + \tau_{jk} + \beta(\theta_k - t_{ik}) - M(1 - x_{ij}) - S_{jk} \\
        t_{jk} &\le t_{ik} + \tau_{jk} + \beta(\theta_k - t_{ik}) + M(1 - x_{ij}) + S_{jk} \\
        \theta_k^{\min} &\le t_{jk} \le \theta_k^{\max} \\
        D_{jk} &\ge t_{jk} - \theta_k \\
        D_{jk} &\ge \theta_k - t_{jk}
        \end{align*}
    \EndFor

    \State Add to route: $\textit{Route} \gets \textit{Route} \cup \{(i,j)\}$

    \State Update cost:
    \[
    \textit{TotalCost} \mathrel{+}= (T_{ij} + \delta_i)x_{ij} + \lambda_1 \sum_k D_{jk} + \lambda_2 \sum_k S_{jk}
    \]

    \State $\textit{Visited} \gets \textit{Visited} \cup \{j\}$,\quad $\textit{Remaining} \gets \textit{Remaining} \setminus \{j\}$
\EndWhile

\Statex \textbf{Phase $n$+1: Return to Warehouse}
\State Let $i \gets$ last visited customer
\State Solve final leg $i \to 0$ and update:
\State $\textit{Route} \gets \textit{Route} \cup \{(i,0)\}$,\quad 
$\textit{TotalCost} \mathrel{+}= (T_{i0} + \delta_i)x_{i0}$

\Statex \textbf{Return:} adaptive route $\{x_{ij}=1\}$, temperatures $\{t_{ik}\}$, and final cost
\end{algorithmic}
\end{algorithm}
\endgroup

\section{Numerical Results}

The foundation of the study is based on real-time monitoring and systematic data collection of environmental conditions affecting perishable produce. Key parameters such as temperature variations, humidity levels, and transit conditions are captured by IoT sensors, which enhances the robustness of the predictive models. Continuous data acquisition is ensured by deploying these sensors across multiple locations, forming the backbone of the optimization models.
\subsection{Scenario Generation}

To systematically evaluate the robustness and practical effectiveness of the proposed adaptive optimization model under real-world uncertainties, we conducted an extensive scenario-based testing framework. The purpose was to simulate diverse operational conditions that fresh food logistics systems routinely encounter, such as random handling delays and ambient temperature fluctuations during transportation. Our approach involved mathematically generating controlled perturbations around key parameters, ensuring that the variability introduced remained statistically justified and representative of actual logistics challenges. Specifically, scenarios were constructed by introducing stochastic variations in two critical parameters: travel delays and ambient temperature shifts. The parameter ranges and probability distributions were selected based on operational data from logistics studies, thermal stability characteristics of common fresh produce, and standard deviation estimates reflective of real-world variability. Normal distributions were chosen due to their well-known properties under the central limit theorem, enabling a realistic aggregation of small environmental factors, whereas truncated Normal-distributions for delay magnitudes ensured non-negativity while preserving natural spread. Each delivery location \( i \in \mathcal{N} \) independently triggered a potential delay occurrence modeled as a Bernoulli random variable:
\[
\mathbb{P}(\text{Delay at stop } i) = p_d,
\]
where \( p_d = 0.2 \) represents a 20\% chance of encountering a delay. If a delay was realized, its magnitude \( \delta_i \) was drawn from a truncated Normal distribution:
\[
\delta_i \sim \max\{0, \mathcal{N}(\mu_d, \sigma_d^2)\},
\]
with mean \( \mu_d = 0.5 \) hours and standard deviation \( \sigma_d = 0.2 \) hours. The truncation at zero ensures that physically infeasible negative delays are avoided, maintaining the realism of operational disturbances. In parallel, ambient temperature shifts were modeled as independent continuous random variables at each delivery stop \( j \) and for each product \( k \in \mathcal{K} \). The temperature change \( \tau_{jk} \) was sampled from a standard Normal-distribution:
\[
\tau_{jk} \sim \mathcal{N}(0, \sigma_\tau^2),
\]
with \( \sigma_\tau = 1.0 \)°C. Fifty independent scenarios were generated under this framework, each representing a distinct realization of travel delays and ambient temperature shifts. For each scenario, the adaptive optimization model was re-solved, and the key outcomes — total travel time, cumulative temperature deviations, and slack usage — were recorded. This repeated solution approach provided empirical insight into the distribution of model performance across random environmental perturbations, allowing us to assess robustness in both average and worst-case senses.\\

\begin{figure}[h]
    \centering
    \fbox{\includegraphics[width=0.6\textwidth]{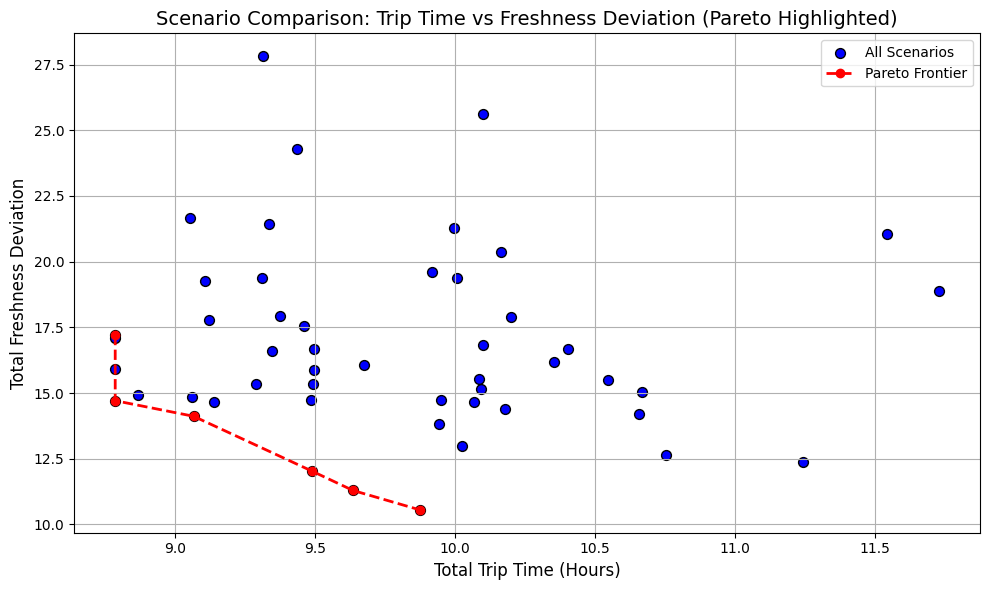}}
    \caption{Scenario-wise comparison of the adaptive optimization model's performance in terms of total trip time versus total freshness deviation.}
    \label{fig:tt}
\end{figure}

We visualized the trade-off between total trip time and cumulative freshness deviation across all tested scenarios. The resulting scatter plot, shown in Figure~\ref{fig:tt}, illustrates the bi-objective performance of the adaptive optimization model, where each blue dot represents a distinct scenario outcome generated from randomized delay and temperature perturbations. The x-axis denotes the total trip time (in hours), while the y-axis indicates the cumulative absolute temperature deviation across all stops and product types.\\

A Pareto frontier is overlaid in red, connecting the set of non-dominated solutions---i.e., those for which no other scenario yields both lower trip time and lower freshness deviation simultaneously. These Pareto-optimal solutions reveal the best-achievable trade-offs under the model’s adaptive control framework. Notably, the frontier is concave and downward-sloping, reflecting the classical tension between rapid delivery and tight temperature control. Scenarios with shorter trip durations (left side of the plot) tend to exhibit higher freshness deviation, suggesting more aggressive routing with less thermal stabilization. In contrast, points near the lower-right corner demonstrate better temperature adherence but at the expense of longer routes, possibly due to conservative pacing or more temperature-conducive paths.\\

This distribution confirms the model’s capacity to balance competing objectives dynamically. The dispersion of the blue points around the frontier also highlights the inherent variability in outcomes under uncertain conditions---validating the need for an adaptive formulation. The frontier itself serves as a decision-support tool, enabling logistics planners to select a delivery policy that aligns with operational priorities (e.g., highlighting freshness over speed in sensitive \glspl{sc}). Furthermore, the tight clustering near the lower-left segment of the frontier suggests that moderate increases in delivery time can yield substantial gains in temperature stability, an insight that could inform risk-aware transportation planning for perishable goods. \\

\begin{figure}[h]
    \centering
    \includegraphics[width=0.8\textwidth]{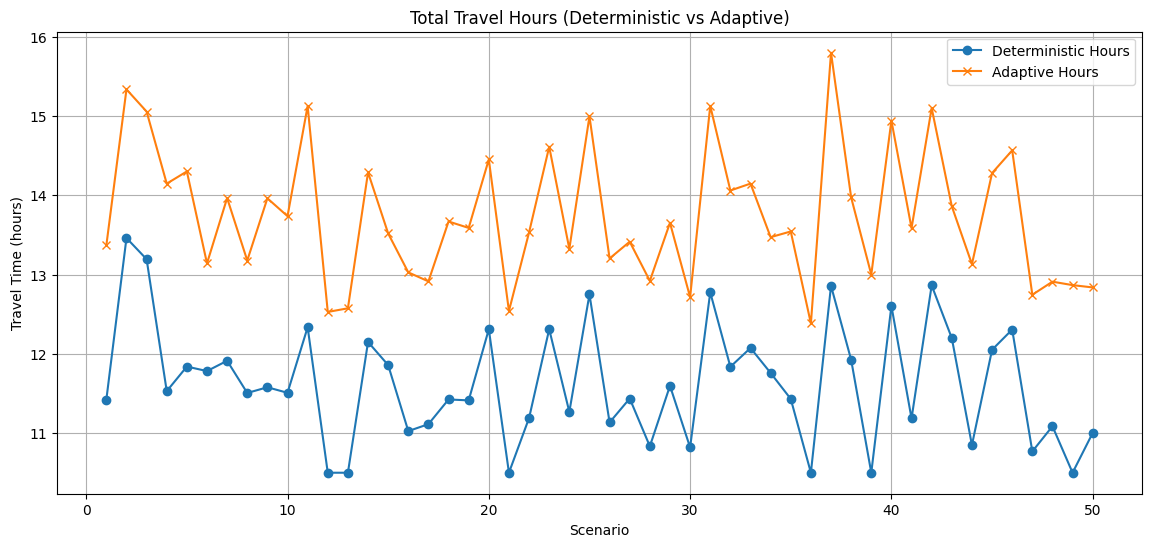}
    \caption{Total Travel Hours (Deterministic vs Adaptive)}
    \label{fig:travel_time_comparison}
\end{figure}

\begin{figure}[h]
    \centering
    \includegraphics[width=0.8\textwidth]{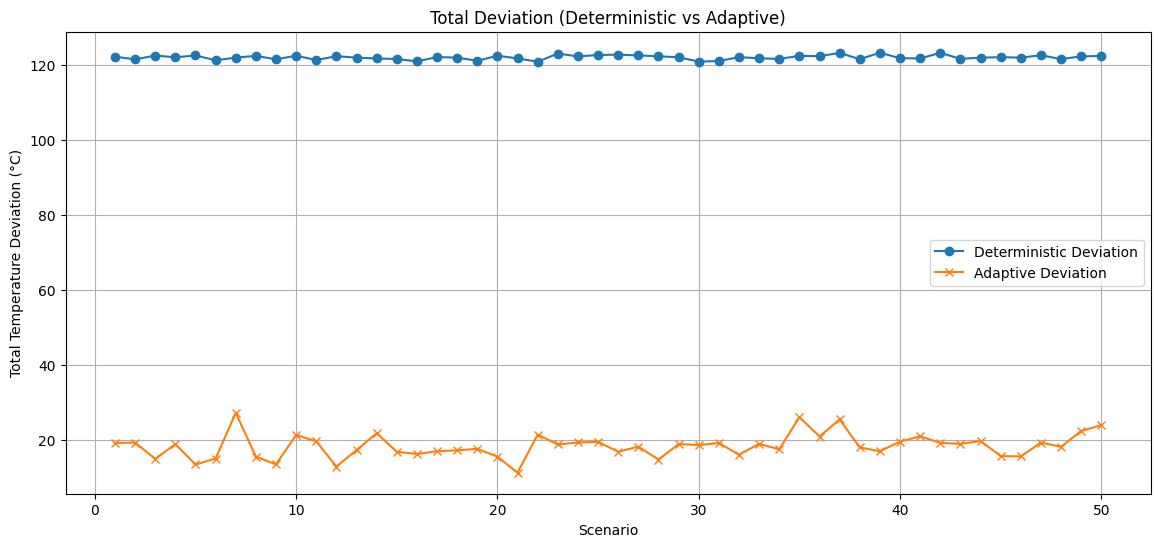}
    \caption{Total Freshness Deviation (Deterministic vs Adaptive)}
    \label{fig:freshness_deviation_comparison}
\end{figure}

We also conducted a scenario-based comparative analysis between a deterministic model (with no freshness penalty) and an adaptive model (which penalizes deviation from ideal product temperature). We executed 50 randomized scenarios using Normal-distributed travel delays and ambient temperature changes, applying both models to each scenario. Two key performance metrics were assessed: total travel time and total freshness deviation. Figure~\ref{fig:travel_time_comparison} presents the travel time (in hours) across 50 scenarios for both models. On the $x$-axis, we plot the scenario index, while the $y$-axis indicates the total route time. The adaptive model consistently incurs slightly higher travel times, ranging approximately from 12.5 to 15.8 hours across scenarios, compared to 10.5 to 13.5 hours in the deterministic model. This increase reflects the routing flexibility required to maintain product temperature closer to ideal values, suggesting a deliberate trade-off for improved freshness.\\
Figure~\ref{fig:freshness_deviation_comparison} illustrates the cumulative temperature deviation from the ideal across all stops and products. The deterministic model exhibits a constant total deviation of approximately 121°C across all scenarios over 10 stops for the four types of produce, reflecting the lack of adaptive correction. In contrast, the adaptive model significantly reduces the deviation (to between 13 and 30 ° C), demonstrating its ability to dynamically regulate internal temperature despite external fluctuations.\\
The deterministic model optimizes solely for travel efficiency, resulting in minimal route time but potentailly large deviations from temperature requirements, this risking spoilage and product degradation. On the other hand, the adaptive model integrates a temperature control mechanism into the optimization, accepting slightly longer travel routes to minimize the temperature deviation penalty. This dual-objective approach is beneficial for logistics involving perishable items, as it ensures compliance with freshness standards while maintaining operational feasibility. Notably, the adaptive model on average achieves up to an 80\% reduction in freshness deviation on average, validating its practical advantage in cold-chain logistics and food supply networks. \\

In addition, to evaluate the performance of the proposed adaptive and deterministic models with respect to product shelf life, we conducted an experiment using a single representative product item—\textit{apple}—which was selected due to its well-documented post-harvest behavior and widespread use in cold chain studies. Apples generally have a practical shelf life of approximately 30 days under ideal cold storage conditions (5\textdegree{}C) \citep{adewoyin2023postharvest}, which we used as the initial product temperature at the start of the delivery route. The modeled delivery trip included 10 stops (including the warehouse), simulating a realistic logistics scenario for fresh produce distribution. During transit, temperature fluctuations within \(\pm 3^\circ\text{C}\) of the ideal were allowed, representing typical cold-chain variation. Shelf life at the end of the route was estimated using a \(Q_{10}\) temperature coefficient of 2, which assumes the degradation rate of the product doubles for every 10\textdegree{}C rise above the optimal temperature. The adaptive model, which responded dynamically to these fluctuations, maintained an average transit temperature of 6.23\textdegree{}C and yielded an estimated shelf life of 27.54 days. In contrast, the deterministic model, with a higher average temperature of 8.68\textdegree{}C, resulted in a lower estimated shelf life of 23.25 days. These results demonstrate the adaptive model’s superior ability to mitigate temperature deviations and preserve product freshness compared to the deterministic baseline.

\subsection{IoT-based Sensors}
To facilitate, Tive Tags, a specialized type of IoT-enabled data logger designed to track temperature fluctuations throughout the transportation process, are utilized \citep{tivetag2025}. These Tags offer an efficient and non-intrusive solution for monitoring perishable goods without requiring complex installation procedures. Tive Tags, which are thin, flexible labels supporting end-to-end temperature monitoring of sensitive products, are suited to the data collection needs due to their ease of use and cost-effectiveness. The process begins with activation, which involves affixing the Tive Tag to the target location or shipment and initiating it by tapping with a Near-field communication (NFC)-enabled smartphone. Once activated, the data logging phase is started, during which temperature data is recorded at regular intervals throughout the monitoring period. Finally, data retrieval is conducted by tapping the Tag again to upload the recorded data to the Tive cloud platform via the mobile application. A comprehensive and secure audit trail of temperature data is ensured by this approach, which is crucial for analyzing the environmental conditions experienced during transit.

\begin{figure}[h]
    \centering
    \fbox{\begin{minipage}{0.45\textwidth}
        \centering
        \includegraphics[width=\textwidth,height=5cm,keepaspectratio]{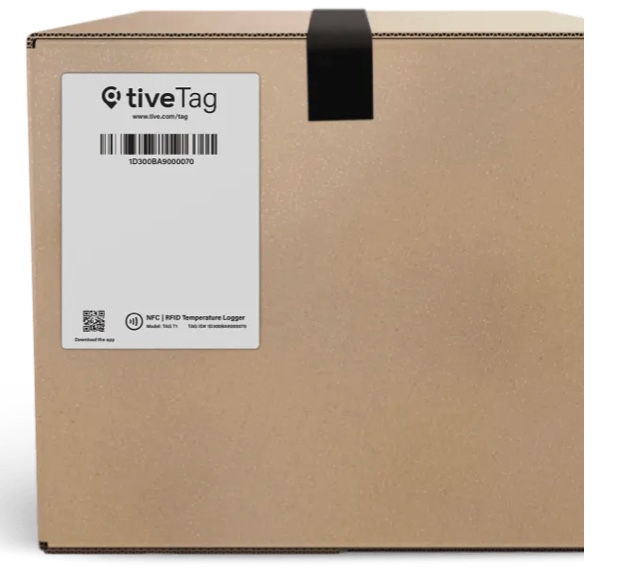}
    \end{minipage}}
    \hspace{5mm}
    \fbox{\begin{minipage}{0.45\textwidth}
        \centering
        \includegraphics[width=\textwidth,height=5cm,keepaspectratio]{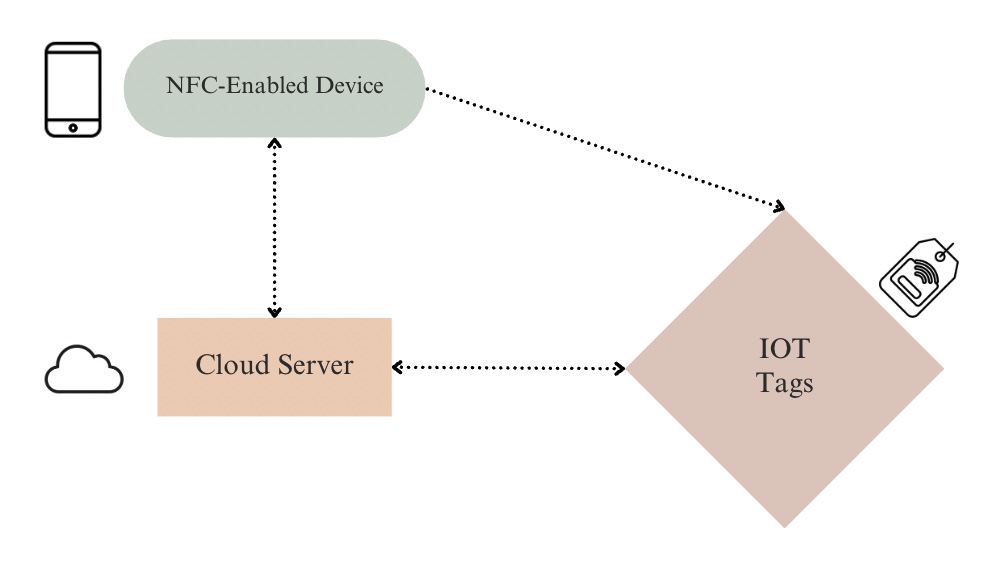}
    \end{minipage}}
    \caption{Left: Tive Tag affixed to a shipment for temperature monitoring. Right: Flowchart illustrating the data movement from the Tive Tag to the cloud and mobile application.}
    \label{fig:tivetagflow}
\end{figure}

\subsection{Synthetic Data Generation}
To evaluate the robustness of the optimization models, synthetic data based on the recorded data from IoT tags that mirror the structural complexity of real-world perishable goods logistics were generated. The distribution network is composed of a central warehouse located at node \(0\) and a set of delivery stops represented by the set \(\mathcal{N} = \{1, 2, \dots, N\}\)


Each delivery stop \(i \in \mathcal{N} \setminus \{0\}\) is assigned a base travel distance \(d_{0i}\) from the warehouse, where
\[
d_{0i} \sim \mathcal{U}(d_{\min}, d_{\max})
\]
for constants \(d_{\min}\) and  \(d_{\max} > 0\). Random delays
\[
\delta_i \sim \mathcal{U}(\delta_{\min}, \delta_{\max})
\]
are introduced at each stop \(i \in \mathcal{N}\) to account for operational disruptions such as traffic congestion or handling issues.

A set of perishable product types denoted by \(\mathcal{K} = \{1, 2, \dots, K\}\) is considered, where each product \(k \in \mathcal{K}\) is characterized by an initial (ideal) shelf life \(L_k^0\). This value denotes the maximum duration (in hours) that the product remains fresh under optimal temperature conditions. The nominal shelf life \(L_k^0\) is perturbed by a random factor\\
\[
\epsilon_k \sim \mathcal{U}(-\alpha, \alpha)
\]
to reflect real-world variability in freshness and spoilage rates:

\begin{equation}
    L_k = L_k^0 (1 + \epsilon_k)
\end{equation}

Furthermore, each product \(k\) is assigned a required shelf life threshold \(R_k\), ensuring the product's viability upon delivery. These constraints are enforced in all models to maintain product integrity.\\

To assess the robustness of different optimization strategies, a consistent synthetic data generation scheme is maintained across all model variants. Additional uncertainties are incorporated in the form of stochastic delays \( \delta_i \), representing variations in traffic, route diversions, or other disruptions. These random variables ensure dynamic travel times \( t_{ij} = d_{ij} / v + \delta_i \), where \( v \) is the nominal vehicle speed. By maintaining controlled randomness in distances, delays, and shelf life parameters, the data framework enables equitable and reproducible comparison across deterministic, robust, DRO, and adaptive models. This allows rigorous evaluation of the impact of real-world uncertainties on delivery cost, route feasibility, and freshness preservation in perishable product distribution.\\

Simulations were conducted based on the synthetic dataset outlined previously. A total of \(K = 4\) perishable produce types—apples, bananas, tomatoes, and strawberries—and a delivery network consisting of \(N = 10\) nodes: one central warehouse located at node \(0\) and \(10\) delivery stops indexed by \(i \in \mathcal{N} \setminus \{0\}\) were considered. Each product \(k \in \mathcal{K}\) was initialized with a nominal shelf life \(L_k^0\) and starting temperature \(\theta_k\), with shelf life perturbed by a uniform distribution factor as defined earlier. For each product, the minimum required shelf life \(R_k\) at the time of delivery is enforced across all model variants.\\

All models were tested using an uniform vehicle speed of \(40\) km/h. The travel times between locations were calculated drawn from a uniform distribution in the range \([15, 60]\) km, and \(\delta_i\) represents the random delay terms described in the data generation section.\\

The results presented below report the total transportation time (in hours) as the primary optimization objective across all model variants. Since each vehicles travels at a constant overall speed, the corresponding total distance traveled can be computed by multiplying the reported objective value by \(40\) km/h.\\

The total trip durations vary depending on the optimization approach used. The same underlying logistics problem is optimized under different treatments of uncertainty by each model. The deterministic model provides a baseline with no uncertainty consideration, while the RO and DRO models introduce increasing levels of conservativeness to handle parameter uncertainty. The adaptive model further incorporates temperature dynamics and active control during transportation, balancing cost-efficiency with produce preservation.

\begin{table}[htbp]
    \centering
    \caption{Total transportation time (in hours) for each model variant}
    \label{tab:model_results}
    \begin{tabular}{lc}
        \toprule
        \textbf{Model} & \textbf{Total Time (hours)} \\
        \midrule
        Deterministic & 9.31 \\
        Robust Optimization (RO) & 15.50 \\
        Stochastic Programming (SP) & 9.42 \\
        Distributionally Robust Optimization (DRO) & 10.50 \\
        Adaptive & 9.80 \\
        \bottomrule
    \end{tabular}
\end{table}

The results presented in Table~\ref{tab:model_results} and the preceding analysis highlight several important findings regarding the performance of different optimization models for perishable product distribution. The deterministic model, serving as a baseline, achieves a competitive total trip time of 9.31 hours but does not account for any uncertainty in travel or environmental conditions. The RO model, designed to handle worst-case parameter uncertainty, results in a significantly longer delivery time of 15.50 hours. This conservativeness ensures reliability but comes at a high cost in terms of efficiency. The SP and DRO models strike a balance between efficiency and uncertainty handling, yielding improved times of 9.42 and 10.50 hours, respectively. The adaptive model achieves overall travel time of 9.80 hours, indicating its effectiveness in dynamically managing environmental and operational factors.\\

One of the key differentiators of the adaptive model is its explicit incorporation of temperature dynamics through continuous monitoring and suitable adjustment. Unlike the RO, DRO, and SP models, which primarily focus on preserving a minimum shelf life via transit time constraints, the adaptive model introduces penalties for temperature deviation and slack compensation directly into the objective function. This enables the model to actively manage freshness throughout the route, rather than passively relying on shelf-life buffers.\\

\subsection{Senstivity Analysis}
To evaluate the performance of the adaptive temperature control model in cold-chain delivery, we also conducted sensitivity analyses on two key parameters: correction strength (\(\beta\)) and ambient temperature shifts (\(\tau_{ik}\)). These govern how actively the system mitigates ambient disturbances and directly influence product freshness across the delivery route.\\

\begin{figure}[h]
    \centering
    \begin{subfigure}{0.48\textwidth}
        \centering
        \includegraphics[width=\linewidth]{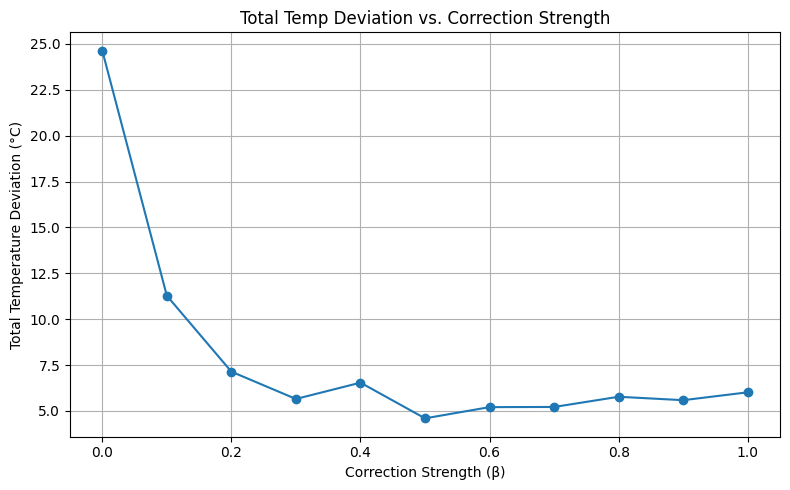}
        \caption{Total temperature deviation vs. correction strength (\(\beta\)).}
        \label{fig:beta_deviation}
    \end{subfigure}
    \hfill
    \begin{subfigure}{0.48\textwidth}
        \centering
        \includegraphics[width=\linewidth]{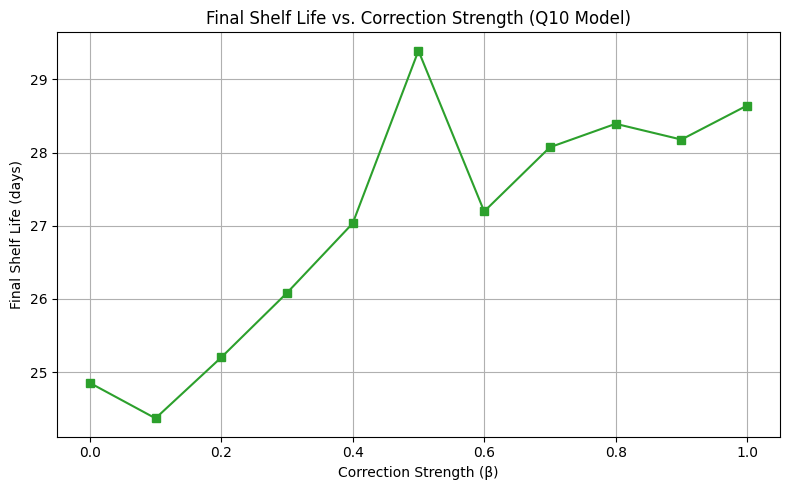}
        \caption{Final shelf life vs. correction strength (\(\beta\)).}
        \label{fig:beta_shelf_life}
    \end{subfigure}
    \caption{Sensitivity analysis of correction strength \(\beta\): (a) Total temperature deviation and (b) final shelf life.}
    \label{fig:beta_sensitivity_grid}
\end{figure}

The model assumes that perishable goods must be maintained near \(5^\circ\mathrm{C}\), within an allowable range of \(2^\circ\mathrm{C}\) to \(8^\circ\mathrm{C}\), with shelf life degradation modeled using the Q$_{10}$ function (\(Q_{10} = 2\)) and a reference shelf life of 30 days. Ambient disturbances are introduced through a uniform distribution \(\mathcal{U}[0, 1]\) for the correction factor \(\beta\), and through normal distribution noise with varying standard deviation \(\sigma_{\tau}\) for ambient temperature shifts \(\tau_{ik}\). Figure~\ref{fig:beta_sensitivity_grid} illustrates the impact of increasing correction strength. As shown in panel~(a), the total temperature deviation decreases sharply with increasing \(\beta\), falling from approximately \(25.5^\circ\mathrm{C}\) at \(\beta = 0.0\) to around \(7.5^\circ\mathrm{C}\) at \(\beta = 0.4\), beyond which it stabilizes around \(5\)–\(6^\circ\mathrm{C}\). This indicates that even modest correction significantly reduces cumulative thermal deviation. The diminishing returns beyond \(\beta = 0.4\) reflect the system’s limited ability to counteract random, bidirectional fluctuations through stronger correction—especially when successive stops experience opposing temperature shifts.

Panel~(b) shows that final shelf life increases with \(\beta\), peaking near 29.5 days around \(\beta = 0.5\), then plateauing. The non-monotonic behavior suggests that overly aggressive correction can lead to minor overcompensations, reducing thermal stability. This aligns with practical cold-chain systems, where energy-efficient controllers or passive cooling systems avoid abrupt changes that may cause oscillations or stress refrigeration units. Thus, a moderate correction strength (\(\beta \in [0.3, 0.5]\)) offers a practical trade-off between energy cost, control effort, and freshness preservation.\\

\begin{figure}[h]
    \centering
    \begin{subfigure}{0.48\textwidth}
        \centering
        \includegraphics[width=\linewidth]{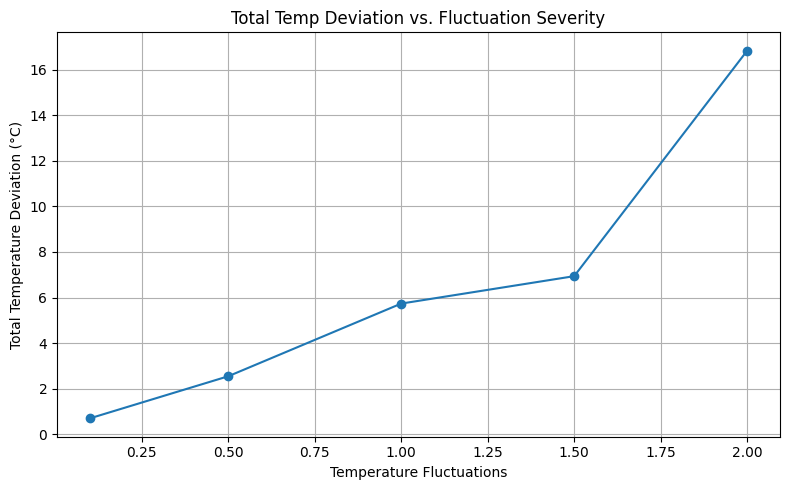}
        \caption{Total temperature deviation vs. fluctuation severity (\(\tau_{ik}\)).}
        \label{fig:temp_deviation}
    \end{subfigure}
    \hfill
    \begin{subfigure}{0.48\textwidth}
        \centering
        \includegraphics[width=\linewidth]{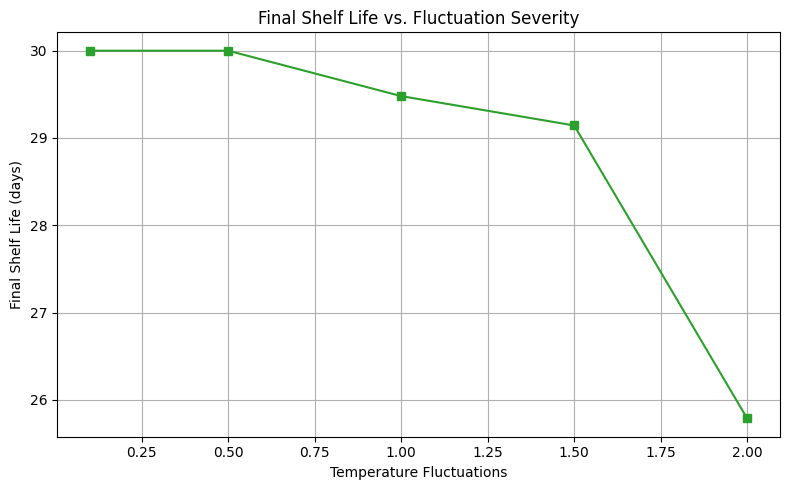}
        \caption{Final shelf life vs. fluctuation severity (\(\tau_{ik}\)).}
        \label{fig:temp_shelf_life}
    \end{subfigure}
    \caption{Sensitivity analysis of fluctuation severity \(\tau_{ik}\): (a) Total temperature deviation and (b) final shelf life.}
    \label{fig:temp_fluctuation_grid}
\end{figure}

Figure~\ref{fig:temp_fluctuation_grid} presents the model's response to increasing ambient temperature shifts (\(\tau_{ik}\)), with \(\beta\) fixed. In panel~(a), the total temperature deviation rises monotonically from ~1\(^\circ\)C at \(\tau_{ik} = 0.25\) to ~16\(^\circ\)C at \(\tau_{ik} = 2.0\), indicating the system’s inability to fully correct large, rapid disturbances even at maximum correction strength. In real-world delivery conditions, this may correspond to frequent door openings, poor insulation, or transit through varying climate zones—factors that cause sharp and compounding thermal stress.\\

Panel~(b) shows the shelf life degradation resulting from these deviations. At low \(\tau_{ik}\), shelf life remains near the reference of 30 days, but it drops sharply at higher \(\tau_{ik}\), reaching ~26 days at \(\tau_{ik} = 2.0\). Given that shelf life decay under the Q$_{10}$ model is exponential in temperature, even small thermal violations have amplified effects. This mirrors practical scenarios where highly sensitive goods (e.g., vaccines, seafood, dairy) require more than algorithmic control—such as insulated containers, real-time temperature tracking, or priority delivery routing to avoid spoilage.\\

Together, these results emphasize that moderate correction strengths are generally sufficient in stable environments, allowing for energy-efficient and effective temperature control. However, in high-volatility scenarios (e.g., long-haul deliveries with varying ambient exposure or cross-border logistics), controller effectiveness saturates. In such cases, passive insulation, proactive rerouting, and IoT-based real-time monitoring become necessary to preserve quality. These findings suggest that operational strategies should jointly optimize controller settings (\(\beta\)) and environmental resilience (mitigating \(\tau_{ik}\)) to ensure freshness and safety in cold-chain logistics.\\

\section{Conclusion}

A dynamic and sensor-aware framework for optimizing the transportation of perishable food products using an adaptive Optimization model is proposed. The model leverages real-time IoT feedback to iteratively refine routing and temperature control decisions across multiple delivery stages. By incorporating shelf-life decay mechanisms grounded in Arrhenius and Q$_{10}$ kinetics, biologically consistent decision-making under time-varying thermal stress is ensured. A comparative evaluation against classical optimization approaches—including deterministic, robust optimization (RO), stochastic programming (SP), and distributionally robust optimization (DRO)—reveals that the adaptive model offers superior performance in both freshness preservation and total travel efficiency. In particular, the model demonstrates an ability to dynamically react to operational uncertainties such as fluctuating ambient temperatures and stop-wise delays, where static formulations fail to capture temporal dependencies and local feedback. Extensive scenario-based testing further validates the model’s reliability under diverse conditions. These simulations, which mimic real-world logistical uncertainties, enable a realistic assessment of the model’s decision robustness. Complementary sensitivity analyses of critical parameters, such as the correction strength \( \beta \) and temperature fluctuation severity, provide valuable insights into the system's behavior across different deployment contexts. These findings reinforce the practical viability of the proposed framework and offer actionable guidelines for system tuning and infrastructure planning in cold-chain logistics. Looking ahead, future extensions will focus on real-time technologies, such as deploying reinforcement learning techniques to support adaptive policy refinement under long-term uncertainty. The generalizability and modularity of the framework also make it well-suited for integration into smart logistics platforms aimed at achieving resilience, sustainability, and freshness assurance across global perishable supply networks.

\backmatter

\bmhead{Acknowledgments}
Not applicable
\bmhead{Declarations}
\begin{itemize}
  \item \textbf{Funding:} Not applicable
  \item \textbf{Conflict of interest:} The authors declare no competing interests.
  \item \textbf{Ethics approval:} Not applicable
  \item \textbf{Consent to participate:} Not applicable
  \item \textbf{Consent for publication:} Not applicable
  \item \textbf{Availability of data:} Available on request
  \item \textbf{Code availability:} Available on request
  \item \textbf{Authors' contributions:} All authors contributed equally.
\end{itemize}

\bibliography{sn-bibliography} 

\end{document}